# HIGHER-ORDER ASYMPTOTIC NORMALITY OF APPROXIMATIONS TO THE MODIFIED SIGNED LIKELIHOOD RATIO STATISTIC FOR REGULAR MODELS[1]

By Heping He and Thomas A. Severini

*Northwestern University*

Approximations to the modified signed likelihood ratio statistic are asymptotically standard normal with error of order $n^{-1}$, where $n$ is the sample size. Proofs of this fact generally require that the sufficient statistic of the model be written as $(\hat{\theta}, a)$, where $\hat{\theta}$ is the maximum likelihood estimator of the parameter $\theta$ of the model and $a$ is an ancillary statistic. This condition is very difficult or impossible to verify for many models. However, calculation of the statistics themselves does not require this condition. The goal of this paper is to provide conditions under which these statistics are asymptotically normally distributed to order $n^{-1}$ without making any assumption about the sufficient statistic of the model.

**1. Introduction.** Let $Y_1, Y_2, \dots, Y_n$ be independent and identically distributed random variables, each distributed according to a distribution with density function $p(y;\theta)$, $\theta \in \Theta$. The likelihood function for $\theta$ based on $Y_1, Y_2, \dots, Y_n$ is given by

$$L(\theta) \equiv L(\theta; Y_1, Y_2, \dots, Y_n) = \prod_{i=1}^{n} p(y_i;\theta), \qquad \theta \in \Theta;$$

let $\ell(\theta) = \log L(\theta)$ denote the log-likelihood function.

Suppose that $\theta = (\psi, \lambda)$, where $\psi$ is a scalar parameter of interest and $\lambda$ is a nuisance parameter. Inference about $\psi$ can be based on the signed likelihood ratio statistic $R$,

$$R \equiv R(\psi) = \operatorname{sgn}(\hat{\psi} - \psi)[2(\ell(\hat{\theta}) - \ell(\hat{\theta}_\psi))]^{1/2},$$

Received July 2005; revised January 2007.
[1]Supported by the NSF under Grants DMS-01-02274 and DMS-06-04123.

where $\hat{\theta}$ denotes the maximum likelihood estimator of $\theta$ and $\hat{\theta}_\psi$ denotes the maximum likelihood estimator of $\theta$ with $\psi$ held fixed. Under standard regularity conditions, $R$ is asymptotically distributed according to a standard normal distribution with error $O(n^{-1/2})$; see, for example, [4].

Adjustments to $R$ have been proposed which reduce this error to $O(n^{-1})$ or even $O(n^{-3/2})$. One such statistic is the modified signed likelihood ratio statistic $R^*$, proposed by [2]. General discussions of the derivation of $R^*$ and its properties are given in [3, 4, 5, 13, 17] and [23]. The analysis in these papers is based on the following assumption:

ASSUMPTION 1.1. The sufficient statistic for the model may be written as $(\hat{\theta}, a)$, where $\hat{\theta}$ denotes the maximum likelihood estimator of $\theta$ and $a$ is an ancillary statistic.

Here ancillarity is taken to mean that the distribution of $a$ does not depend on $\theta$, although if $a$ is used in Assumption 1.1, it must satisfy the additional requirement that $(\hat{\theta}, a)$ is sufficient. Under Assumption 1.1, the log-likelihood function $\ell(\theta)$ may be written as $\ell(\theta; \hat{\theta}, a)$ to emphasize that it is a function of the sufficient statistic $(\hat{\theta}, a)$ as well as the usual argument $\theta$. Note that some distributions, such as full-rank exponential family distributions and certain transformation models, satisfy this assumption; however, for other types of models, it is often very difficult or impossible to verify Assumption 1.1.

The modified signed likelihood ratio statistic $R^*$ is given by

$$R^* = R + \frac{1}{R}\log\left(\frac{U}{R}\right),$$

where

$$U = \left|\begin{matrix} (\ell_{\cdot;\hat{\theta}}(\hat{\theta}) - \ell_{\cdot;\hat{\theta}}(\hat{\theta}_\psi))^T \\ \ell_{\lambda;\hat{\theta}}(\hat{\theta}_\psi) \end{matrix}\right| \Big/ (|j_{\lambda\lambda}(\hat{\theta}_\psi)|^{1/2}|\hat{J}|^{1/2}),$$

$$\ell_{\cdot;\hat{\theta}}(\theta) = \frac{\partial \ell(\theta; \hat{\theta}, a)}{\partial \hat{\theta}}, \qquad \ell_{\lambda;\hat{\theta}}(\hat{\theta}_\psi) = \frac{\partial^2 \ell(\hat{\theta}_\psi; \hat{\theta}, a)}{\partial \lambda\, \partial \hat{\theta}},$$

$$j_{\lambda\lambda}(\hat{\theta}_\psi) = -\ell_{\lambda\lambda}(\hat{\theta}_\psi) \quad \text{and} \quad \hat{J} = -\ell_{\theta\theta}(\hat{\theta}; \hat{\theta}, a).$$

The sign of $U$ is taken to be the same as that of $R$ so that $U/R$ is never negative. In the moderate-deviation case, $R^*$ is asymptotically normally distributed with error $O(n^{-3/2})$, while in the large-deviation case $R^*$ is asymptotically normally distributed with error $O(n^{-1})$; see [2]. Note that the errors in both cases are uniform and relative.

For models in which Assumption 1.1 cannot be verified, computation of $U$ is difficult or impossible. [10] uses a variable $T$ in place of $U$; Assumption 1.1 is not needed for determination of $T$. However, it does require a prior density function for $\theta$. [22] and [18] propose, respectively, an approximation $\bar{R}^*$ to $R^*$ based on covariances and an approximation $\hat{R}^*$ to $R^*$ based on empirical covariances.

The covariance-based approximation is based on the fact that, in the expression for $U$, $\ell_{\cdot;\hat{\theta}}(\hat{\theta}) - \ell_{\cdot;\hat{\theta}}(\theta)$ and $\ell_{\theta;\hat{\theta}}(\theta)$ may be approximated by

$$\bar{\ell}_{\cdot;\hat{\theta}}(\hat{\theta}) - \bar{\ell}_{\cdot;\hat{\theta}}(\theta) = \{Q(\hat{\theta};\hat{\theta}) - Q(\theta;\hat{\theta})\} i(\hat{\theta})^{-1}\hat{J}$$

and $\bar{\ell}_{\theta;\hat{\theta}}(\theta) = I(\theta;\hat{\theta}) i(\hat{\theta})^{-1}\hat{J}$ respectively, where

$$I(\theta;\theta_\circ) = E[\ell_\theta(\theta)\ell_\theta(\theta_\circ)^T;\theta_\circ], \qquad Q(\theta;\theta_\circ) = E[\ell(\theta)\ell_\theta^T(\theta_\circ);\theta_\circ]$$

and $i(\hat{\theta}) = I(\hat{\theta};\hat{\theta})$. Let $\bar{U}$ denote the approximation of the statistic $U$ based on the above approximation; then an approximation to $R^*$ is given by

$$\bar{R}^* = R + \frac{1}{R}\log\frac{\bar{U}}{R};$$

see [22] for further details.

Let

$$\hat{Q}(\theta;\theta_\circ) = \sum \ell^{(j)}(\theta)\ell_\theta^{(j)}(\theta_\circ)^T, \qquad \hat{I}(\theta;\theta_\circ) = \sum \ell_\theta^{(j)}(\theta)\ell_\theta^{(j)}(\theta_\circ)^T$$

and $\hat{\imath}(\hat{\theta}) = \hat{I}(\hat{\theta};\hat{\theta})$, where $\ell^{(j)}(\theta)$ denotes the log-likelihood function based on the $j$th observation alone. Then $\ell_{\cdot;\hat{\theta}}(\hat{\theta}) - \ell_{\cdot;\hat{\theta}}(\theta)$ and $\ell_{\theta;\hat{\theta}}(\theta)$ may be approximated by

$$\hat{\ell}_{\cdot;\hat{\theta}}(\hat{\theta}) - \hat{\ell}_{\cdot;\hat{\theta}}(\theta) = \{\hat{Q}(\hat{\theta};\hat{\theta}) - \hat{Q}(\theta;\hat{\theta})\}\hat{\imath}(\hat{\theta})^{-1}\hat{J}$$

and $\hat{\ell}_{\theta;\hat{\theta}}(\theta) = \hat{I}(\theta;\hat{\theta})\hat{\imath}(\hat{\theta})^{-1}\hat{J}$ respectively. Let $\hat{U}$ denote the statistic $U$ based on these approximations; then an approximation to $R^*$ is given by

$$\hat{R}^* = R + \frac{1}{R}\log(\hat{U}/R);$$

see [18] for further details.

The higher-order asymptotic normality of $\bar{R}^*$ and $\hat{R}^*$ has been established by [22] and [18], respectively, under Assumption 1.1. The goal of this paper is to prove that these statistics are asymptotically normal with error $O(n^{-1})$ without requiring Assumption 1.1.

## 2. Notation and regularity conditions.

For $n = 1, 2, \ldots,$ let $Y = (Y_1, Y_2, \ldots, Y_n)$, where $Y_1, Y_2, \ldots, Y_n$ are i.i.d. random vectors with each of dimension $d_\circ$ and let $p(\cdot, \theta)$ denote the density function of $Y_1$. Here $\theta$ is a $k$-dimensional parameter of the form

$$\theta = (\theta_1, \theta_2, \ldots, \theta_k) \equiv (\psi, \lambda_1, \lambda_2, \ldots, \lambda_{k-1}) \equiv (\psi, \lambda),$$

where $\psi$ is the scalar parameter of interest and $\lambda$ is the nuisance parameter vector. Let $L(\theta) \equiv L(\theta; Y)$ denote the likelihood function for $\theta \in \Theta$ and $Y \in \mathcal{R}^{n \times d_\circ}$ and let $\ell(\theta) = \log L(\theta)$. For each $\theta \in \Theta$, let

$$E_\theta = \{Y \in \mathcal{R}^{n \times d_\circ} | L(\theta; Y) > 0\}.$$

For $r, s = 1, 2, \ldots, k$, define

$$\ell_\theta(\theta) = \frac{\partial}{\partial \theta} \ell(\theta), \qquad \ell_{\theta\theta}(\theta) = \frac{\partial^2}{\partial \theta \, \partial \theta^T} \ell(\theta),$$

$$\ell_{\theta\theta\theta_r}(\theta) = \frac{\partial^3}{\partial \theta \, \partial \theta^T \, \partial \theta_r} \ell(\theta) \quad \text{and} \quad \ell_{\theta\theta\theta_r\theta_s}(\theta) = \frac{\partial^4}{\partial \theta \, \partial \theta^T \, \partial \theta_r \, \partial \theta_s} \ell(\theta).$$

Assumption 2.1. Assume that $\Theta$ is an open subset of $\mathcal{R}^k$.

1. $\ell(\theta)$ is five times continuously differentiable at each $\theta \in \Theta$ with respect to $\theta_i$, $i = 1, 2, \ldots, k$, for all $Y \in E_\theta$;
2. $\int \frac{\partial}{\partial \theta_i} L(\theta; y) \, dy = \int \frac{\partial^2}{\partial \theta_i \, \partial \theta_j} L(\theta; y) \, dy = \int \frac{\partial^3}{\partial \theta_i \, \partial \theta \, \partial \theta_\kappa} L(\theta; y) \, dy = \int \frac{\partial^4}{\partial \theta_i \, \partial \theta_j \, \partial \theta_\kappa \, \partial \theta_\ell} L(\theta; y) dy = 0; i, j, \kappa, \ell = 1, 2, \ldots, k$, for each $\theta \in \Theta$;
3. $\mathrm{Var}\{\ell_\psi(\theta)\} = E[\ell_\psi^2(\theta)] > 0$ , $\mathrm{Var}\{\ell_\lambda(\theta)\} = E[\ell_\lambda(\theta) \ell_\lambda(\theta)^T]$ is positive definite, for each $\theta \in \Theta$;
4. $(\ell_{\theta_i}(\theta), \ell_{\theta_j\theta_k}(\theta), \ell_{\theta_\ell\theta_m\theta_p}(\theta), \ell_{\theta_q\theta_r\theta_s\theta_t}(\theta))$ has joint cumulants up to third order, where $i, j, \kappa, \ell, m, p, q, r, s, t = 1, 2, \ldots, k$;
5. The maximum likelihood estimator of $\theta$ is consistent as $n \to \infty$.

Conditions ensuring that part (5) of Assumption 2.1 holds are given, for example, by [1, 24] and [16].

Let $\ell^{(i)}(\theta) = \log p(Y_i; \theta)$, $i = 1, 2, \ldots, n$. Then, for $r, s, t = 1, \ldots, k$, we have

$$\ell(\theta) = \sum_{i=1}^n \ell^{(i)}(\theta), \qquad \ell_\theta(\theta) = \sum_{i=1}^n \ell_\theta^{(i)}(\theta),$$

$$\ell_{\theta\theta}(\theta) = \sum_{i=1}^n \ell_{\theta\theta}^{(i)}(\theta), \qquad \ell_{\theta_r\theta_s\theta_t}(\theta) = \sum_{i=1}^n \ell_{\theta_r\theta_s\theta_t}^{(i)}(\theta).$$

For $r = 1, \ldots, k$, let

$$Z_\theta = \sum_{i=1}^n \frac{1}{\sqrt{n}} \ell_\theta^{(i)}(\theta),$$

$$Z_{\theta\theta} = \sum_{i=1}^n \frac{1}{\sqrt{n}} (\ell_{\theta\theta}^{(i)}(\theta) - V_{\theta\theta}),$$

$$Z_{\theta\theta\theta_r} = \sum_{i=1}^n \frac{1}{\sqrt{n}} (\ell_{\theta\theta\theta_r}^{(i)}(\theta) - V_{\theta\theta\theta_r}),$$

where $V_{\theta\theta} = E[\ell^{(i)}_{\theta\theta}(\theta)], V_{\theta\theta\theta_r} = E[\ell^{(i)}_{\theta\theta\theta_r}(\theta)]$, $i = 1, 2, \dots, n$. Let

$$\begin{aligned} Z = (&Z_\theta, Z_{\psi\theta}, Z_{\lambda_1\lambda}, Z_{\lambda_2\lambda}, \dots, Z_{\lambda_{k-1}\lambda}, Z_{\psi\psi\psi}, Z_{\psi\psi\lambda}, \\ &Z_{\psi\lambda_1\lambda}, Z_{\psi\lambda_2\lambda}, \dots, Z_{\psi\lambda_{k-1}\lambda}, Z_{\lambda_1\lambda_1\lambda}, Z_{\lambda_1\lambda_2\lambda}, \dots, \\ &Z_{\lambda_1\lambda_{k-1}\lambda}, Z_{\lambda_2\lambda_2\lambda}, Z_{\lambda_2\lambda_3\lambda}, \dots, Z_{\lambda_2\lambda_{k-1}\lambda}, \dots, Z_{\lambda_{k-1}\lambda_{k-1}\lambda})^T. \end{aligned}$$

Let $d$ denote the rank of the covariance matrix of $Z$. By the spectral decomposition of the covariance matrix of $Z$, there exist a matrix $T = (t_{ij})_{d\times(k+k^2+k^3)}$ and a random vector $X$ with $d$ elements such that $X = T \cdot Z$, $\mathrm{Cov}(X) = I_{d\times d}$ and $Z = (T'T)^{-1}T'X$. Let

$$\tilde{Y}_i = \begin{pmatrix} t_{1,1}\ell^{(i)}_\psi(\theta) + t_{1,2}\ell^{(i)}_{\lambda_1}(\theta) + \\ \quad \cdots + t_{1,k+k^2+k^3}(\ell^{(i)}_{\lambda_{k-1}\lambda_{k-1}\lambda_{k-1}}(\theta) - V_{\lambda_{k-1}\lambda_{k-1}\lambda_{k-1}}) \\ t_{2,1}\ell^{(i)}_\psi(\theta) + t_{2,2}\ell^{(i)}_{\lambda_1}(\theta) + \\ \quad \cdots + t_{2,k+k^2+k^3}(\ell^{(i)}_{\lambda_{k-1}\lambda_{k-1}\lambda_{k-1}}(\theta) - V_{\lambda_{k-1}\lambda_{k-1}\lambda_{k-1}}) \\ \vdots \\ t_{d,1}\ell^{(i)}_\psi(\theta) + t_{d,2}\ell^{(i)}_{\lambda_1}(\theta) + \\ \quad \cdots + t_{d,k+k^2+k^3}(\ell^{(i)}_{\lambda_{k-1}\lambda_{k-1}\lambda_{k-1}}(\theta) - V_{\lambda_{k-1}\lambda_{k-1}\lambda_{k-1}}) \end{pmatrix},$$

where $i = 1, 2, \dots, n$. Obviously, we have $X = \sum_{i=1}^n \frac{1}{\sqrt{n}}\tilde{Y}_i$, $\mathrm{Cov}(\tilde{Y}_i) = I_{d\times d}$.

Let $\|\cdot\|$ denote the Euclidean norm.

ASSUMPTION 2.2. The following three points are assumed:

I. The distribution of $\tilde{Y}_i$ has a finite cumulant generating function in a neighborhood of the origin, that is, $|\log(E[\exp\{<t, \tilde{Y}_i>\}])| < +\infty$, where $\|t\| < \varepsilon$ for some $\varepsilon > 0$.

II. The characteristic function $\xi$ of $\tilde{Y}_i$ satisfies Cramér's condition, that is, for some $k > 0$, $\sup\{|\xi(t)|; \|t\| > k\} < 1$.

III. Some power of $\xi$ is integrable, that is, $\int_{-\infty}^{+\infty} \xi(t)^\alpha\, dt < +\infty$ for some real number $\alpha > 1$.

Condition II can be verified by applying the Riemann–Lebesgue lemma; see, for example, [7]. Condition III can be verified by applying a result such as Lemma 2.4.3 in [14].

## 3. Main result.

The following theorem is the main result of the paper.

THEOREM 3.1. *Assume that the Assumptions 2.1 and 2.2 hold. Then, uniformly for $-\infty < t < \infty$, we have*

$$Pr(\bar{R}^* \le t) = \int_{-\infty}^t \frac{1}{\sqrt{2\pi}} \exp\left(-\frac{x^2}{2}\right) dx\, [1 + O(n^{-1})]$$

*and*

$$Pr(\hat{R}^* \leq t) = \int_{-\infty}^{t} \frac{1}{\sqrt{2\pi}} \exp\left(-\frac{x^2}{2}\right) dx\,[1 + O(n^{-1})].$$

*That is, $\bar{R}^*$ and $\hat{R}^*$ converge in distribution to standard normal random variables with the error of order of $n^{-1}$.*

Note that the $O(n^{-1})$ errors in the normal approximations to the distribution functions of $\bar{R}^*$ and $\hat{R}^*$ are relative and uniform in $t$, the argument of the distribution function under consideration. The uniformity follows directly from Theorem 3.2 of [19], while relativeness follows from the fact that the cumulative distribution function of a standard normal random variable is $O(1)$ as $n \to \infty$.

The proof of Theorem 3.1 is based on the following lemma.

LEMMA 3.1. *Assume that Assumption 2.2 is satisfied. Define $X$ as in Section 2; then*

$$Pr(X \in \mathcal{B}) = \int_{\mathcal{B}} \xi_n(u)\, du + o(n^{-1}),$$

*where $E\{\|X\|^4\} < +\infty$; $\mathcal{B}$ is any Borel set in $\mathcal{R}^d$; $\xi_n(u) = \sum_{r=0}^{2} P_r(-\phi : \{\chi_{v,n}\})(u)$; $P_r(-\phi : \{\chi_{v,n}\})(u)$ is the density function of the finite signed measure with characteristic function $\tilde{P}_r(it : \{\chi_{v,n}\}) \exp\{-\frac{1}{2}\|t\|^2\}$, $t \in R^d$, $r = 0, 1, 2$; $\tilde{P}_r(z : \{\chi_{v,n}\})$, $r = 0, 1, 2$, are Cramér–Edgeworth polynomials defined in [6]; and $\chi_{v,n}, 1 \leq v \leq 4$ are the cumulants of $X$, satisfying $\chi_{1,n} = 0, \chi_{2,n} = 1_{\mathcal{R}^d}$. That is, equations (3.1)–(3.3) in [19] are satisfied for $s = 4$.*

PROOF. [21] gives conditions for the validity of multivariate Edgeworth expansions for a sequence of random variables. We check those conditions so that we can use [21] to prove this lemma.

Using the properties of $X$ given in the previous section, condition I of Assumption 2.2 implies condition (i) in Corollary 4.2 of [21] and condition II of Assumption 2.2 implies condition (ii) in Corollary 4.2 of [21]. It is trivial to verify that condition (iii) in Corollary 4.2 of [21] is also satisfied. Condition III in Assumption 2.2 implies that some power of $\xi$ is integrable. Therefore, the conditions of Theorem 3.1 in [21] for $s = 4, \rho_{4,n}(t) = a_n(t)^{-1} = n^{-1/2}, \varepsilon_n = n^{-3/2}$ follow from Corollary 4.2 in [21].

Applying Theorem 3.1 in [21] shows that equation (3.7) in Theorem 3.1 in [21] holds for $s = 4, \varepsilon_n = n^{-3/2}$ and, hence, that condition in Corollary (3.3) of [21] is satisfied. Therefore, equation (3.10) in Corollary (3.3) of [21] gives that

$$Pr(X \in \mathcal{B}) = \int_{\mathcal{B}} \xi_n(u)\, du + o(n^{-1}),$$

where $\mathcal{B}$ and $\xi_n(u)$ are the notation explained in this lemma. By Assumption 2.1, it is obvious that $E\{\|X\|^4\} < \infty$. Therefore, equations (3.1), (3.2) and (3.3) in [19] hold for $s = 4$. $\square$

Proof of Theorem 3.1. We begin by outlining the basic idea of the proof. Suppose that a sequence of distributions can be approximated to a certain order of accuracy by an Edgeworth series; then Theorem 3.2 of Skovgaard [19] shows that such an expansion may be transformed by a sequence of deterministic or stochastic smooth functions of the corresponding random vectors to yield another Edgeworth expansion of the resulting sequence of distributions. The terms in this expansion may be calculated from the moments obtained by the delta method, that is, the moments formally calculated from a Taylor series expansion omitting terms of high order. Thus, to prove the theorem, we use Lemma 3.1 to obtain an Edgeworth series approximation to the distribution of $X$ and then, in Section 5, show that $\bar{\bar{R}}^*$ has an expansion in terms of $X$. Theorem 3.2 of [19] then gives an Edgeworth series approximation to the distribution of $\bar{\bar{R}}^*$; the same approach can be used for $\hat{R}^*$. As noted above, the uniformity in Theorem 3.1 follows from the uniformity in equation (3.13) of Theorem 3.2 of [19].

Lemma 3.1 shows that equations (3.1)–(3.3) in [19] hold. Note that $\bar{R}^*$ does not satisfy the first part of equation (3.4) of [19]; thus, we begin by constructing a location adjustment of $\bar{R}^*$. Let

$$\bar{\bar{R}}^* \equiv \bar{R}^* + [\tfrac{1}{6}(-V_{\psi\psi})^{-3/2} V_{\psi\psi\psi} - \tfrac{1}{2}(-V_{\psi\psi})^{-1/2} \operatorname{tr}(V_{\lambda\lambda}^{-1} V_{\lambda\lambda\psi}) + \tfrac{1}{2}(-V_{\psi\psi})^{-3/2} v_{\psi,\psi\psi} - (-V_{\psi\psi})^{-1/2} \operatorname{tr}(V_{\lambda\lambda}^{-1} v_{\lambda,\lambda\psi})] n^{-1/2}.$$

Expanding the corresponding likelihood equations shows that $\hat{\theta}$ and $\hat{\theta}_\psi$ can be considered as stochastic functions of $Z$; see [20], particularly equation (6.2). The formula of $\bar{\bar{R}}^*$ contains $\hat{\theta}$ and $\hat{\theta}_\psi$. By its expansion, $\bar{\bar{R}}^*$ can be similarly considered as a stochastic function of $Z$ and it is equal to 0 when $Z = 0$. Hence, $\bar{\bar{R}}^*$ satisfies the first part of equation (3.4) in [19]. However, $\bar{R}^*$ considered as a stochastic function of $Z$ is not equal to zero when $Z = 0$. Hence, Theorem 3.2 in [19] can be applied only to $\bar{\bar{R}}^*$ to obtain its asymptotic distribution. Note that the function $f_n$ in Theorem 3.2 in [19] can be a stochastic function; see Remark 3.8 in [19]. The asymptotic distribution of $\bar{R}^*$ is obtained by reversing the location adjustment in the final part of this proof.

By using Taylor series expansions and collecting terms by orders of $n^{-1/2}$, $\bar{\bar{R}}^*$ may be expanded,

$$\bar{\bar{R}}^* = (-V_{\psi\psi})^{-1/2} Z_\psi + [\tfrac{1}{2}(-V_{\psi\psi})^{-3/2} Z_{\psi\psi} Z_\psi - (-V_{\psi\psi})^{-1/2} Z_{\psi\lambda}^T V_{\lambda\lambda}^{-1} Z_\lambda$$

$$(3.1)\qquad - \tfrac{1}{2}(-V_{\psi\psi})^{-3/2}V^T_{\psi\psi\lambda}V^{-1}_{\lambda\lambda}Z_\lambda Z_\psi + \tfrac{1}{2}(-V_{\psi\psi})^{-1/2}Z^T_\lambda V^{-1}_{\lambda\lambda}V_{\psi\lambda\lambda}V^{-1}_{\lambda\lambda}Z_\lambda + \tfrac{1}{6}(-V_{\psi\psi})^{-5/2}V_{\psi\psi\psi}Z^2_\psi]n^{-1/2} + R(n^{-1}) + O_p(n^{-3/2}),$$

where $R(n^{-1})$ represents a finite number of terms of order $n^{-1}$ such that the first four cumulants of $nR(n^{-1})$ are bounded in $n$; see Section 5 for details of the derivation of equation (3.1). Even though we can obtain the specific form of $R(n^{-1})$ by expanding $\bar{R}^*$ to order $n^{3/2}$, we omit its specific form since it is extremely tedious. Note that $Z_\psi, Z_\lambda, \ldots$ in equation (3.1) may be viewed as functions of $X$. Thus, we can write $\bar{\bar{R}}^* = f(X,\omega)$, where $f(\cdot;\omega)$ represents the function corresponding to the first three terms in equation (3.1); $\omega$ denotes the sample point of the underlying experiment and accounts for the contribution of the $O_p$ term in (1). It follows that Theorem 3.2 in [19] can be applied to $\bar{\bar{R}}^*$; see Remark 3.8 of [19].

We know from equation (3.1) that at $X = 0$, $\bar{\bar{R}}^*$ is four times differentiable, $\bar{\bar{R}}^*(0) = 0$ and $D\bar{\bar{R}}^*(0)$ is of rank one; it follows that equation (3.4) in [19] holds. It is shown below that $\bar{\bar{R}}^*$ has asymptotic variance 1; hence, $\bar{\bar{R}}^*$ does not have to be normalized as in equation (3.5) of [19].

Write $\bar{\bar{R}}^* = \bar{\bar{R}} + O_p(n^{-3/2})$, where $\bar{\bar{R}}$ denotes the first three terms on the right-hand side of equation (3.1). Note that $\bar{\bar{R}}$ is a function of $X$. By the uniqueness of power series expansions, $\bar{\bar{R}}$ agrees with the Taylor series expansion of $\bar{\bar{R}}^*$ with respect to $X$ up to a polynomial of a finite number of terms, neglecting terms of order $O_p(n^{-3/2})$.

Using the expression for $\bar{\bar{R}}$ given in equation (3.1), we obtain the following expression for $E[\bar{\bar{R}}]$:

$$\begin{aligned} E[\bar{\bar{R}}] = \{&\tfrac{1}{2}(-V_{\psi\psi})^{-3/2}v_{\psi,\psi\psi} \\ &- (-V_{\psi\psi})^{-1/2}E[Z^T_{\psi\lambda}V^{-1}_{\lambda\lambda}Z_\lambda] \\ &+ \tfrac{1}{2}(-V_{\psi\psi})^{-3/2}V^T_{\psi\psi\lambda}V^{-1}_{\lambda\lambda}V_{\psi\lambda} \\ &+ \tfrac{1}{2}(-V_{\psi\psi})^{-1/2}E[Z^T_\lambda V^{-1}_{\lambda\lambda}V_{\psi\lambda\lambda}V^{-1}_{\lambda\lambda}Z_\lambda] \\ &- \tfrac{1}{6}(-V_{\psi\psi})^{-5/2}V_{\psi\psi\psi}V_{\psi\psi}\}n^{-1/2} + O(n^{-1}) \\ = \{&\tfrac{1}{2}(-V_{\psi\psi})^{-3/2}v_{\psi,\psi\psi} + \tfrac{1}{6}(-V_{\psi\psi})^{-3/2}V_{\psi\psi\psi} \\ &- (-V_{\psi\psi})^{-1/2}\operatorname{tr}(V^{-1}_{\lambda\lambda}v_{\lambda,\lambda\psi}) \\ &- \tfrac{1}{2}(-V_{\psi\psi})^{-1/2}\operatorname{tr}(V^{-1}_{\lambda\lambda}V_{\lambda\lambda\psi})\}n^{-1/2} + O(n^{-1}). \end{aligned}$$

Similar computations show that

$$V(\bar{\bar{R}}) = 1 + O(n^{-1}),$$

$$K_3(\bar{\bar{R}}) = O(n^{-1}) \quad \text{and} \quad K_4(\bar{\bar{R}}) = O(n^{-1}).$$

Let $\eta_n$ be the density function with characteristic function

$$\hat{\eta}_n(t) = \exp\{i\langle t, \tilde{\chi}_{1,n}\rangle - \tfrac{1}{2}\tilde{\chi}_{2,n}(t,t)\} \sum_{r=0}^{2} \tilde{P}_r(it : \{\tilde{\chi}_{v,n}\}),$$

where

$$\begin{aligned}\tilde{\chi}_{1,n} = [&\tfrac{1}{2}(-V_{\psi\psi})^{-3/2} v_{\psi,\psi\psi} + \tfrac{1}{6}(-V_{\psi\psi})^{-3/2} V_{\psi\psi\psi} \\ &- (-V_{\psi\psi})^{-1/2}\operatorname{tr}(V_{\lambda\lambda}^{-1} v_{\lambda,\lambda\psi}) \\ &\qquad - \tfrac{1}{2}(-V_{\psi\psi})^{-1/2}\operatorname{tr}(V_{\lambda\lambda}^{-1} V_{\lambda\lambda\psi})] n^{-1/2} + O(n^{-1}), \\ \tilde{\chi}_{2,n} = 1 &+ O(n^{-1}), \qquad \tilde{\chi}_{3,n} = O(n^{-1}), \qquad \tilde{\chi}_{4,n} = O(n^{-1}).\end{aligned}$$

Define

$$\rho_n(\alpha) = (\log n)^{1/2}(2+\alpha)^{1/2}, \qquad \alpha > 0,$$

and

$$H_n(\alpha) = \{t \in \mathcal{R}^{k+k^2+k^3} \mid \|t\| \le \rho_n(\alpha)\}, \qquad \alpha > 0.$$

We now verify Assumption 3.1 in [19]. For all $\alpha > 0$ and all sufficiently large $n$, we have the following:

1. $\bar{\bar{R}}^*$ is four times continuously differentiable on $H_n(\alpha)$, and
$$\sup\{\|D^4 \bar{\bar{R}}^*(t)\| \mid t \in H_n(\alpha)\} = o(n^{-1}).$$
2. Define $\lambda_n = \sup\{(\|D^j \bar{\bar{R}}^*(0)\|/j!)^{1/(j-1)} | 2 \le j \le 3\}$, then $\lambda_n^3 = o(n^{-1})$.

The continuous differentiability of $\bar{\bar{R}}^*$ claimed in point 1 requires some clarification. Note that $\bar{R}^*$ has a singularity at points where $R = 0$ and $\bar{U} = 0$. However, the singularity is removable. Whenever $R = 0$, we must have $\psi = \hat{\psi}$, from which $\bar{U} = 0$ follows. When $\psi$ approaches $\hat{\psi}$, equations (5.11) and (5.16) in Section 5 show that $R = O((\hat{\psi} - \psi))$ and equations (5.7), (5.8) and (5.11) show that $\bar{U} = O((\hat{\psi} - \psi))$. Hence, $\bar{U}$ and $R$ approach zero at the same rate whenever $R$ goes to zero. Furthermore, equation (5.14) shows that $\bar{U}/R$ goes to 1. Therefore, the singularity is removed.

Therefore, we have shown that all conditions in Theorem 3.2 in [19] are satisfied for $q = p = s = 4$; it follows that

$$\eta_n(t) = \left[\sum_{r=0}^{2} P_r(-\phi_{0,\tilde{\chi}_{2,n}} : \{\tilde{\chi}_{v,n}\})(x)\right]\Big|_{x = x - \tilde{\chi}_{1,n}},$$

where

$$P_0(-\phi_{0,\tilde{\chi}_{2,n}} : \{\tilde{\chi}_{v,n}\})(x) = \frac{1}{\sqrt{2\pi}\sqrt{\tilde{\chi}_{2,n}}} \exp\Big(-\frac{x^2}{2\tilde{\chi}_{2,n}}\Big),$$

$$P_1(-\phi_{0,\tilde{\chi}_{2,n}} : \{\tilde{\chi}_{v,n}\})(x) = -\frac{\tilde{\chi}_{3,n}}{6}\frac{1}{\sqrt{2\pi}\sqrt{\tilde{\chi}_{2,n}}}\Big[\frac{3x}{\tilde{\chi}_{2,n}^2}\exp\Big(-\frac{x^2}{2\tilde{\chi}_{2,n}}\Big) - \frac{x^3}{\tilde{\chi}_{2,n}^3}\exp\Big(-\frac{x^2}{2\tilde{\chi}_{2,n}}\Big)\Big],$$

$$\begin{aligned} P_2(-\phi_{0,\tilde{\chi}_{2,n}} : \{\tilde{\chi}_{v,n}\})(x) &= \frac{\tilde{\chi}_{4,n}}{24\sqrt{2\pi}\sqrt{\tilde{\chi}_{2,n}}}\Big[\frac{3}{\tilde{\chi}_{2,n}^2}\exp\Big(-\frac{x^2}{2\tilde{\chi}_{2,n}}\Big) - \frac{6x^2}{\tilde{\chi}_{2,n}^3}\exp\Big(-\frac{x^2}{2\tilde{\chi}_{2,n}}\Big) \\ &\quad + \frac{x^4}{\tilde{\chi}_{2,n}^4}\exp\Big(-\frac{x^2}{2\tilde{\chi}_{2,n}}\Big)\Big] \\ &\quad + \frac{\tilde{\chi}_{3,n}^2}{72\sqrt{2\pi}\sqrt{\tilde{\chi}_{2,n}}}\Big[-\frac{15}{\tilde{\chi}_{2,n}^3}\exp\Big(-\frac{x^2}{2\tilde{\chi}_{2,n}}\Big) + \frac{45x^2}{\tilde{\chi}_{2,n}^4}\exp\Big(-\frac{x^2}{2\tilde{\chi}_{2,n}}\Big) \\ &\quad - \frac{15x^4}{\tilde{\chi}_{2,n}^5}\exp\Big(-\frac{x^2}{2\tilde{\chi}_{2,n}}\Big) + \frac{x^6}{\tilde{\chi}_{2,n}^6}\exp\Big(-\frac{x^2}{2\tilde{\chi}_{2,n}}\Big)\Big]. \end{aligned}$$

Let

$$\begin{aligned} C \equiv [&\tfrac{1}{6}(-V_{\psi\psi})^{-3/2}V_{\psi\psi\psi} - \tfrac{1}{2}(-V_{\psi\psi})^{-1/2}\operatorname{tr}(V_{\lambda\lambda}^{-1}V_{\lambda\lambda\psi}) \\ &+ \tfrac{1}{2}(-V_{\psi\psi})^{-3/2}v_{\psi,\psi\psi} - (-V_{\psi\psi})^{-1/2}\operatorname{tr}(V_{\lambda\lambda}^{-1}v_{\lambda,\lambda\psi})]n^{-1/2}. \end{aligned}$$

Finally, it follows from equation (3.13) in [19] that

$$\begin{aligned} Pr(\bar{R}^* \le t) &= Pr(\bar{\bar{R}}^* \le t + C) \\ &= \int_{-\infty}^{C+t}\sum_{r=0}^{2} P_r(-\phi_{0,\tilde{\chi}_{2,n}} : \{\tilde{\chi}_{v,n}\})|_{x=x-C}\,dx + o\Big(\frac{1}{n}\Big) \\ &= \int_{-\infty}^{t}\frac{1}{\sqrt{2\pi}\sqrt{\tilde{\chi}_{2,n}}}\exp\Big(-\frac{x^2}{2\tilde{\chi}_{2,n}}\Big)\,dx \\ &\quad + \tilde{\chi}_{3,n}\int_{-\infty}^{t} -\frac{1}{6}\frac{1}{\sqrt{2\pi}\sqrt{\tilde{\chi}_{2,n}}}\Big[\frac{3x}{\tilde{\chi}_{2,n}^2}\exp\Big(-\frac{x^2}{2\tilde{\chi}_{2,n}}\Big) \end{aligned}$$

$$
\begin{aligned}
&\qquad\qquad - \frac{x^3}{\tilde{\chi}_{2,n}^3} \exp\Big(-\frac{x^2}{2\tilde{\chi}_{2,n}}\Big)\Big] \, dx \\
&+ \tilde{\chi}_{4,n} \int_{-\infty}^{t} \frac{1}{24\sqrt{2\pi}\sqrt{\tilde{\chi}_{2,n}}} \Big[\frac{3}{\tilde{\chi}_{2,n}^2} \exp\Big(-\frac{x^2}{2\tilde{\chi}_{2,n}}\Big) \\
&\qquad - \frac{6x^2}{\tilde{\chi}_{2,n}^3} \exp\Big(-\frac{x^2}{2\tilde{\chi}_{2,n}}\Big) \\
&\qquad + \frac{x^4}{\tilde{\chi}_{2,n}^n} \exp\Big(-\frac{x^2}{2\tilde{\chi}_{2,n}}\Big)\Big] \, dx \\
&+ \tilde{\chi}_{3,n}^2 \int_{-\infty}^{t} \frac{1}{72\sqrt{2\pi}\sqrt{\tilde{\chi}_{2,n}}} \Big[-\frac{15}{\tilde{\chi}_{2,n}^3} \exp\Big(-\frac{x^2}{2\tilde{\chi}_{2,n}}\Big) \\
&\qquad + \frac{45x^2}{\tilde{\chi}_{2,n}^4} \exp\Big(-\frac{x^2}{2\tilde{\chi}_{2,n}}\Big) \\
&\qquad - \frac{15x^4}{\tilde{\chi}_{2,n}^5} \exp\Big(-\frac{x^2}{2\tilde{\chi}_{2,n}}\Big) \\
&\qquad + \frac{x^6}{\tilde{\chi}_{2,n}^6} \exp\Big(-\frac{x^2}{2\tilde{\chi}_{2,n}}\Big)\Big] \, dx + o\Big(\frac{1}{n}\Big).
\end{aligned}
$$

Note that

$$
\tilde{\chi}_{2,n} = 1 + O\Big(\frac{1}{n}\Big), \qquad \tilde{\chi}_{3,n} = O\Big(\frac{1}{n}\Big), \qquad \tilde{\chi}_{4,n} = O\Big(\frac{1}{n}\Big).
$$

Therefore, we have

$$
Pr(\bar{R}^* \le t) = \int_{-\infty}^{t} \frac{1}{2\pi} \exp\Big(-\frac{x^2}{2}\Big) \, dx + O\Big(\frac{1}{n}\Big).
$$

The proof for $\hat{R}^*$ is very similar. We can expand $\hat{R}^*$ to the third order and the $k + k^2 + k^3$ elements in the vector $Z$ appear in its expansion as well. Since $\hat{R}^* = \bar{R}^* + O_p(1/n)$, the differences between the expansions of $\hat{R}^*$ and $\bar{R}^*$ exist only in terms of order $n^{-1}$ or higher. The expansion of $\hat{R}^*$ has the first, second, third and fourth order cumulants equal to the ones of the expansion of $\bar{R}^*$, up to the order $O(n^{-1})$. Therefore, similar computations give

$$
Pr(\hat{R}^* \le t) = \int_{-\infty}^{t} \frac{1}{2\pi} \exp\Big(-\frac{x^2}{2}\Big) \, dx + O\Big(\frac{1}{n}\Big),
$$

which finally concludes this proof. □

TABLE 1
*Remission times for patients*

| | | | | | | | | | | |
|---|---|---|---|---|---|---|---|---|---|---|
| 1 | 1 | 2 | 2 | 3 | 4 | 4 | 5 | 5 | 6 | 8 |
| 8 | 9 | 10 | 10 | 12 | 14 | 16 | 20 | 24 | 34 | |

**4. An example.** The linear-exponential model is used in the analysis of survival data. Under this model, survival times are assumed to have a distribution with density of the form

$$f(y) = (\psi + \lambda y) \exp[-(\psi y + \tfrac{1}{2}\lambda y^2)], \qquad y > 0,$$

where $\psi > 0$ and $\lambda > 0$ are unknown parameters; note that this is not an exponential family model. The log-likelihood function, based on the independent survival times $y_1, \ldots, y_n$, is given by

$$\ell(\theta) = \sum_{i=1}^{n} \log(\psi + \lambda y_i) - \sum_{i=1}^{n} (\psi y_i + \tfrac{1}{2}\lambda y_i^2).$$

For this model, it is very difficult or impossible to verify Assumption 1.1. However, Assumptions 2.1 and 2.2 can be quite easily checked. Consider Assumption 2.1. Part 1 of this assumption is naturally satisfied. Part 2 follows from applying Corollary 2.4.1 in [8] on interchanging the order of differentiation and integration. Direct computation can be used to establish part 3; since the cumulant-generating function of the first four log-likelihood derivatives exists, part 4 of Assumption 2.1 is satisfied as well. Now consider Assumption 2.2. Since $\tilde{Y}_i$ has a cumulant-generating function for this model, part I is satisfied. Applying the Riemann–Lebesgue lemma in [7] and Lemma 2.4.3 in [14] establishes parts II and III, respectively.

Thus, it follows from Theorem 2.1 that the modified versions of the signed likelihood ratio statistic, $\bar{R}^*$ and $\hat{R}^*$, are each asymptotically normally distributed with error of order $o(n^{-1})$.

The data in Table 1 are from Example 7.2 in [15]; the values represent remission times in weeks for 21 patients with acute leukemia. Using this data and the linear-exponential model, confidence limits for $\psi$ based on $R$, $\bar{R}^*$ and $\hat{R}^*$ were calculated; these are presented in Table 2. Table 2 also contains the results of a small simulation study designed to estimate the coverage probabilities of confidence limits based on $R$, $\bar{R}^*$ and $\hat{R}^*$. In the simulation study, the data were drawn independently from a linear-exponential distribution with the parameter values taken to be the maximum likelihood estimates based on the data in Table 1. A Monte Carlo sample size of 100,000 was used.

**5. Asymptotic expansion of $\bar{R}^*$.** We now consider the expansion $\bar{R}^*$ used in the proof of Theorem 3.1. Recall that, for $\theta = (\theta_1, \theta_2, \ldots, \theta_k)$, $\psi = \theta_1$ and $\lambda = (\theta_2, \theta_3, \ldots, \theta_k)$. Without loss of generality, we assume that $\psi$ and $\lambda$ are orthogonal parameters; see, for example, [9, 11, 12] and [4] (Section 2.7).

We begin by expanding $(\bar{\ell}_{\cdot;\hat{\theta}}(\hat{\theta}) - \bar{\ell}_{\cdot;\hat{\theta}}(\hat{\theta}_\psi))^T$:

$$(\bar{\ell}_{\cdot;\hat{\theta}}(\hat{\theta}) - \bar{\ell}_{\cdot;\hat{\theta}}(\hat{\theta}_\psi))^T = (\hat{\theta} - \hat{\theta}_\psi)^T \hat{J} - \tfrac{1}{2} \sum_{s,t=1}^{k} \hat{\delta}_s \hat{\delta}_t \mu^T_{\theta,\theta_s\theta_t}(\hat{\theta}) i^{-1}(\hat{\theta}) \hat{J} + O_p(n^{-1/2}), \tag{5.1}$$

where $\mu_{\theta,\theta_s\theta_t}(\theta) = E(\ell_\theta(\theta)\ell_{\theta_s\theta_t}(\theta); \theta)$, $\hat{\delta} = \hat{\theta} - \hat{\theta}_\psi$. The approximation $\bar{\ell}_{\theta;\hat{\theta}}(\hat{\theta}_\psi)$ to $\ell_{\theta;\hat{\theta}}(\hat{\theta}_\psi)$ can be expanded,

$$\bar{\ell}_{\theta;\hat{\theta}}(\hat{\theta}_\psi) = \hat{J} - \sum_{t=1}^{k} \mu_{\theta,\theta\theta_t}(\hat{\theta}) \hat{\delta}_t i^{-1}(\hat{\theta}) \hat{J} + O_p(1), \tag{5.2}$$

so that $\bar{\ell}_{\lambda;\hat{\theta}}(\hat{\theta}_\psi)$ can be expanded,

$$\bar{\ell}_{\lambda;\hat{\theta}}(\hat{\theta}_\psi) = \hat{J}_{\lambda\theta} - \sum_{t=1}^{k} \hat{\delta}_t \mu_{\lambda,\theta\theta_t}(\hat{\theta}) i^{-1}(\hat{\theta}) \hat{J} + O_p(1). \tag{5.3}$$

Now let us turn to Taylor series expansion of the determinant of the asymmetric matrix $A + X$, where $A = (a_{ij})$ and $X = (x_{ij})$ are two square matrices of the same order. Then

$$|A + X| = |A| + \sum_{i,j} \frac{\partial |A + X|}{\partial x_{ij}}\Big|_{x=0} x_{ij} + O(x_{st}^2)$$

TABLE 2

*Upper confidence limits based on $R$, $\bar{R}^*$ and $\hat{R}^*$*

| | Confidence limit | | | Coverage probability | | |
|---|---|---|---|---|---|---|
| **Probability** | $R$ | $\bar{R}^*$ | $\hat{R}^*$ | $R$ | $\bar{R}^*$ | $\hat{R}^*$ |
| 0.010 | 0.0206 | 0.0260 | 0.0263 | 0.0032 | 0.0125 | 0.0126 |
| 0.025 | 0.0293 | 0.0353 | 0.0356 | 0.0150 | 0.0261 | 0.0264 |
| 0.050 | 0.0373 | 0.0439 | 0.0443 | 0.0397 | 0.0526 | 0.0526 |
| 0.100 | 0.0472 | 0.0545 | 0.0550 | 0.0745 | 0.1100 | 0.1102 |
| 0.900 | 0.1327 | 0.1372 | 0.1375 | 0.8840 | 0.8990 | 0.8991 |
| 0.950 | 0.1441 | 0.1469 | 0.1472 | 0.9366 | 0.9467 | 0.9469 |
| 0.975 | 0.1540 | 0.1557 | 0.1559 | 0.9661 | 0.9742 | 0.9748 |
| 0.990 | 0.1657 | 0.1664 | 0.1667 | 0.9808 | 0.9900 | 0.9902 |

$$
\begin{aligned}
&= |A| + \sum_{i,j} \frac{\partial |A+X|}{\partial (a_{ij}+x_{ij})}\Big|_{x=0} x_{ij} + O(x_{st}^2) \\
&= |A| + \sum_{i,j} (-1)^{i+j} |(A+X)_{ij}|\Big|_{x=0} x_{ij} + O(x_{st}^2) \qquad (5.4) \\
&= |A| + \sum_{i,j} (-1)^{i+j} |A_{ij}| x_{ij} + O(x_{st}^2) \\
&= |A| + \operatorname{tr}(\operatorname{adj}(A)X) + O(x_{st}^2) \qquad \text{for any } s, t,
\end{aligned}
$$

where $|(A+X)_{ij}|$ is the minor of $a_{ij}+x_{ij}$ in $|A+X|$; $|A_{ij}|$ is the minor of $a_{ij}$ in $|A|$; $\operatorname{adj}(A)$ is the adjugate matrix of $A$; $\operatorname{tr}(\operatorname{adj}(A)X)$ is the trace of $\operatorname{adj}(A)X$.

Recall that $\bar{U}$ is given by

$$
\bar{U} = \begin{vmatrix} (\bar{\ell}_{\cdot;\hat{\theta}}(\hat{\theta}) - \bar{\ell}_{\cdot;\hat{\theta}}(\hat{\theta}_\psi))^T \\ \bar{\ell}_{\lambda;\hat{\theta}}(\hat{\theta}_\psi) \end{vmatrix} |j_{\lambda\lambda}(\hat{\theta}_\psi)|^{-1/2} |\hat{J}|^{-1/2}. \tag{5.5}
$$

Substituting equations (5.1) and (5.3) into equation (5.5) and applying equation (5.4) gives

$$
\begin{aligned}
\bar{U} = &\begin{vmatrix} (\hat{\theta} - \hat{\theta}_\psi)^T \hat{J} - \frac{1}{2} \sum_{s,t=1}^{k} \hat{\delta}_s \hat{\delta}_t \mu^T_{\theta,\theta_s\theta_t}(\hat{\theta}) i^{-1}(\hat{\theta}) \hat{J} \\ \hat{J}_{\lambda\theta} - \sum_{t=1}^{k} \hat{\delta}_t \mu_{\lambda,\theta\theta_t}(\hat{\theta}) i^{-1}(\hat{\theta}) \hat{J} \end{vmatrix} \\
&\times |j_{\lambda\lambda}(\hat{\theta}_\psi)|^{-1/2} |\hat{J}|^{-1/2} + O_p(n^{-1}).
\end{aligned} \tag{5.6}
$$

Let

$$
A = \begin{pmatrix} (\hat{\theta} - \hat{\theta}_\psi)^T \hat{J} \\ \hat{J}_{\lambda\theta} \end{pmatrix}, \qquad X = \begin{pmatrix} -\frac{1}{2} \sum_{s,t=1}^{k} \hat{\delta}_s \hat{\delta}_t \mu^T_{\theta,\theta_s\theta_t}(\hat{\theta}) i^{-1}(\hat{\theta}) \hat{J} \\ - \sum_{t=1}^{k} \hat{\delta}_t \mu_{\lambda,\theta\theta_t}(\hat{\theta}) i^{-1}(\hat{\theta}) \hat{J} \end{pmatrix} \tag{5.7}
$$

and note that $A+X$ is asymmetric. Applying equation (5.4) once more produces

$$
\begin{aligned}
\bar{U} &= |A+X| |j_{\lambda\lambda}(\hat{\theta}_\psi)|^{-1/2} |\hat{J}|^{-1/2} + O_p(n^{-1}) \\
&= \begin{vmatrix} (\hat{\theta} - \hat{\theta}_\psi)^T \hat{J} \\ \hat{J}_{\lambda\theta} \end{vmatrix} |j_{\lambda\lambda}(\hat{\theta}_\psi)|^{-1/2} |\hat{J}|^{-1/2} \\
&\quad + \operatorname{tr}(\operatorname{adj}(A)X) |j_{\lambda\lambda}(\hat{\theta}_\psi)|^{-1/2} |\hat{J}|^{-1/2} + O_p(n^{-1}).
\end{aligned} \tag{5.8}
$$

Let

$$B \equiv \begin{vmatrix} (\hat{\theta}-\hat{\theta}_\psi)^T \hat{J} \\ \hat{J}_{\lambda\theta} \end{vmatrix} |j_{\lambda\lambda}(\hat{\theta}_\psi)|^{-1/2}|\hat{J}|^{-1/2}R^{-1}, \tag{5.9}$$

which is the first term of equation (5.8) divided by $R$. Using the formula for the determinant of a partitioned matrix, an alternative form for $B$ is

$$\begin{aligned} B &= |\hat{J}_{\lambda\lambda}|[(\hat{\theta}-\hat{\theta}_\psi)^T\hat{J}_{\theta\psi} - (\hat{\theta}-\hat{\theta}_\psi)^T\hat{J}_{\theta\lambda}\hat{J}_{\lambda\lambda}^{-1}\hat{J}_{\lambda\psi}] \\ &\quad \times |j_{\lambda\lambda}(\hat{\theta}_\psi)|^{-1/2}|\hat{J}_{\lambda\lambda}|^{-1/2}R^{-1}[\hat{J}_{\psi\psi} - \hat{J}_{\psi\lambda}\hat{J}_{\lambda\lambda}^{-1}\hat{J}_{\lambda\psi}]^{-1/2} \\ &= \frac{|\hat{J}_{\lambda\lambda}|^{1/2}(\hat{\psi}-\psi)[\hat{J}_{\psi\psi} - \hat{J}_{\psi\lambda}\hat{J}_{\lambda\lambda}^{-1}\hat{J}_{\lambda\psi}]^{1/2}}{R|j_{\lambda\lambda}(\hat{\theta}_\psi)|^{1/2}}. \end{aligned} \tag{5.10}$$

Hence, in order to find expansion of $B$, we have to find an expansion of $R^{-1}$. Note that

$$\hat{\lambda} - \hat{\lambda}_\psi = -(\hat{\psi}-\psi)\hat{J}_{\lambda\lambda}^{-1}\hat{J}_{\lambda\psi} - \tfrac{1}{2}(\hat{\psi}-\psi)^2\hat{J}_{\lambda\lambda}^{-1}\ell_{\lambda\psi\psi}(\hat{\theta}) + O_p(n^{-3/2}). \tag{5.11}$$

Therefore, we have

$$\begin{aligned} \frac{1}{R} &= \operatorname{sgn}(\hat{\psi}-\psi)\Bigg\{|\hat{\psi}-\psi|^{-1}\hat{J}_{\psi\psi}^{-1/2} \\ &\qquad + \frac{1}{6}|\hat{\psi}-\psi|^{-3}\hat{J}_{\psi\psi}^{-3/2} \\ &\qquad \times \left[-\sum_{r,s,t=1}^{k}\ell_{\theta_r\theta_s\theta_t}(\hat{\theta})\hat{\delta}_r\hat{\delta}_s\hat{\delta}_t\right] + O_p(n^{-1})\Bigg\} \\ &= (\hat{\psi}-\psi)^{-1}\hat{J}_{\psi\psi}^{-1/2} \\ &\quad + \frac{1}{6}(\hat{\psi}-\psi)^{-3}\hat{J}_{\psi\psi}^{-3/2} \\ &\quad \times \left[-\sum_{r,s,t=1}^{k}\ell_{\theta_r\theta_s\theta_t}(\hat{\theta})\hat{\delta}_r\hat{\delta}_s\hat{\delta}_t\right] + O_p(n^{-1}). \end{aligned} \tag{5.12}$$

Using

$$\begin{aligned} |j_{\lambda\lambda}(\hat{\theta}_\psi)|^{-1/2} &= |\hat{J}_{\lambda\lambda}|^{-1/2} - \tfrac{1}{2}|\hat{J}_{\lambda\lambda}|^{-3/2}\operatorname{tr}\left(\operatorname{adj}(\hat{J}_{\lambda\lambda})\left(\sum_{s=1}^{k}\ell_{\lambda\lambda\theta_s}(\hat{\theta})\hat{\delta}_s\right)\right) \\ &\qquad + O_p(n^{-(k/2)-(1/2)}), \end{aligned} \tag{5.13}$$

together with equations (5.12), we obtain the following expansion of $B$:

$$B = \frac{1}{R}\Bigg\{(\hat{\psi}-\psi)[\hat{J}_{\psi\psi} - \hat{J}_{\psi\lambda}\hat{J}_{\lambda\lambda}^{-1}\hat{J}_{\lambda\psi}]^{1/2}$$

$$
\begin{aligned}
&- \frac{1}{2}|\hat{J}_{\lambda\lambda}|^{-1}(\hat{\psi}-\psi)[\hat{J}_{\psi\psi}-\hat{J}^T_{\psi\lambda}\hat{J}^{-1}_{\lambda\lambda}\hat{J}_{\lambda\psi}]^{1/2} \\
&\qquad \times \operatorname{tr}\left(\operatorname{adj}(\hat{J}_{\lambda\lambda})\left(\sum_{s=1}^{k}\ell_{\lambda\lambda\theta_s}(\hat{\theta})\hat{\delta}_s\right)\right)+O_p(n^{-1})\Bigg\} \qquad (5.14)\\
&= 1-\frac{1}{2}|\hat{J}_{\lambda\lambda}|^{-1}\operatorname{tr}\left(\operatorname{adj}(\hat{J}_{\lambda\lambda})\left(\sum_{s=1}^{k}\ell_{\lambda\lambda\theta_s}(\hat{\theta})\hat{\delta}_s\right)\right)\\
&\quad +\frac{1}{6}(\hat{\psi}-\psi)^{-2}\hat{J}^{-1}_{\psi\psi}\left[-\sum_{r,s,t=1}^{k}\ell_{\theta_r\theta_s\theta_t}(\hat{\theta})\hat{\delta}_r\hat{\delta}_s\hat{\delta}_t\right]+O_p(n^{-1}).
\end{aligned}
$$

Now consider expansion of $\bar{R}^*$. Note that

$$
\begin{aligned}
\bar{R}^* &= R+\frac{1}{R}\log\frac{\bar{U}}{R}\\
&= R+\frac{1}{R}\log\left\{B+\frac{\operatorname{tr}(\operatorname{adj}(A)X)}{R|j_{\lambda\lambda}(\hat{\theta}_\psi)|^{1/2}|\hat{J}|^{1/2}}+O_p(n^{-1})\right\}\\
&= R-\frac{1}{2}(\hat{\psi}-\psi)^{-1}|\hat{J}_{\lambda\lambda}|^{-1}\hat{J}^{-1/2}_{\psi\psi}\operatorname{tr}\left(\operatorname{adj}(\hat{J}_{\lambda\lambda})\left(\sum_{r=1}^{k}\ell_{\lambda\lambda_r}(\hat{\theta})\hat{\delta}^r\right)\right) \qquad (5.15)\\
&\quad +\frac{1}{6}(\hat{\psi}-\psi)^{-3}\hat{J}^{-3/2}_{\psi\psi}\left[-\sum_{r,s,t=1}^{k}\ell_{\theta_r\theta_s\theta_t}(\hat{\theta})\hat{\delta}_r\hat{\delta}_s\hat{\delta}_t\right]\\
&\quad +\frac{1}{R}\frac{\operatorname{tr}(\operatorname{adj}(A)X)}{\begin{vmatrix}(\hat{\theta}-\hat{\theta}_\psi)^T\hat{J}\\ \hat{J}_{\lambda\theta}\end{vmatrix}}+O_p(n^{-1}).
\end{aligned}
$$

Using the expansion

$$
\begin{aligned}
R &= \operatorname{sgn}(\hat{\psi}-\psi)\Bigg\{[(\hat{\theta}-\hat{\theta}_\psi)^T\hat{J}(\hat{\theta}-\hat{\theta}_\psi)]^{1/2}\\
&\qquad +\tfrac{1}{6}[(\hat{\theta}-\hat{\theta}_\psi)^T\hat{J}(\hat{\theta}-\hat{\theta}_\psi)]^{-1/2}\\
&\qquad \times\left[\sum_{r,s,t=1}^{k}\ell_{\theta_r\theta_s\theta_t}(\hat{\theta})\hat{\delta}_r\hat{\delta}_s\hat{\delta}_t\right]\Bigg\}+O_p(n^{-1}) \qquad (5.16)\\
&= (\hat{\psi}-\psi)\hat{J}^{1/2}_{\psi\psi}+\tfrac{1}{6}(\hat{\psi}-\psi)^{-1}\hat{J}^{-1/2}_{\psi\psi}\\
&\qquad \times\sum_{r,s,t=1}^{k}\ell_{\theta_r\theta_s\theta_t}(\hat{\theta})\hat{\delta}_r\hat{\delta}_s\hat{\delta}_t+O_p(n^{-1})
\end{aligned}
$$

for $R$ in equation (5.15) yields

$$
\begin{aligned}
\bar{R}^* = {} & (\hat{\psi}-\psi)\hat{J}_{\psi\psi}^{1/2} \\
& + \Biggl\{ \frac{1}{6}(\hat{\psi}-\psi)^{-1}\hat{J}_{\psi\psi}^{-1/2} \sum_{r,s,t=1}^{k} \ell_{\theta_r\theta_s\theta_t}(\hat{\theta})\hat{\delta}_r\hat{\delta}_s\hat{\delta}_t \\
& \quad - \frac{1}{6}(\hat{\psi}-\psi)^{-3}\hat{J}_{\psi\psi}^{-3/2} \sum_{r,s,t=1}^{k} \ell_{\theta_r\theta_s\theta_t}(\hat{\theta})\hat{\delta}_r\hat{\delta}_s\hat{\delta}_t \\
& \quad - \frac{1}{2}(\hat{\psi}-\psi)^{-1}|\hat{J}_{\lambda\lambda}|^{-1}\hat{J}_{\psi\psi}^{-1/2} \\
& \quad \times \operatorname{tr}\left(\operatorname{adj}(\hat{J}_{\lambda\lambda})\left(\sum_{r=1}^{k}\ell_{\lambda\lambda_r}(\hat{\theta})\hat{\delta}^r\right)\right) \\
& \quad + \frac{1}{R}\frac{\operatorname{tr}(\operatorname{adj}(A)X)}{\begin{vmatrix} (\hat{\theta}-\hat{\theta}_\psi)^T\hat{J} \\ \hat{J}_{\lambda\theta}\end{vmatrix}} \Biggr\} + O_p(n^{-1}).
\end{aligned}
\tag{5.17}
$$

In order to obtain an expansion of $\bar{R}^*$ in terms of log-likelihood derivatives, each term in equation (5.17) can be expanded. Let

$$
\begin{aligned}
Z_\theta &\equiv \frac{1}{\sqrt{n}}\ell_\theta(\theta), \\
Z_{\theta\theta} &\equiv \frac{1}{\sqrt{n}}(\ell_{\theta\theta}(\theta) - nV_{\theta\theta}), \qquad V_{\theta\theta} = \frac{1}{n}E[\ell_{\theta\theta}(\theta);\theta].
\end{aligned}
\tag{5.18}
$$

Note that we have

$$
Z_{\psi\lambda} = \frac{1}{\sqrt{n}}\ell_{\psi\lambda}(\theta) \tag{5.19}
$$

because of orthogonal parameters. We continue to set

$$
Z_{\theta_r\theta_s\theta_t} = \frac{1}{\sqrt{n}}(\ell_{\theta_r\theta_s\theta_t}(\theta) - nV_{\theta_r\theta_s\theta_t}), \tag{5.20}
$$

where $V_{\theta_r\theta_s\theta_t} = \frac{1}{n}E[\ell_{\theta_r\theta_s\theta_t}(\theta);\theta]$, $r,s,t = 1,2,\dots,k$.

The expansion of the first term of equation (5.17) of $\bar{R}^*$ is given by

$$
\begin{aligned}
&(\hat{\psi}-\psi)\hat{J}_{\psi\psi}^{1/2} \\
&\quad = (-V_{\psi\psi})^{1/2}Z_\psi \\
&\qquad + \biggl[\frac{1}{2}(-V_{\psi\psi})^{-3/2}Z_{\psi\psi}Z_\psi
\end{aligned}
$$

$$
\begin{aligned}
&- (-V_{\psi\psi})^{-1/2} Z_{\psi\lambda}^T V_{\lambda\lambda}^{-1} Z_\lambda \\
&- \frac{1}{2}(-V_{\psi\psi})^{-3/2} V_{\psi\psi\lambda}^T V_{\lambda\lambda}^{-1} Z_\lambda Z_\psi \\
&+ \frac{1}{2}(-V_{\psi\psi})^{-1/2} Z_\lambda^T V_{\lambda\lambda}^{-1} V_{\psi\lambda\lambda} V_{\lambda\lambda}^{-1} Z_\lambda \Big] \frac{1}{\sqrt{n}} \\
&+ O_p(n^{-1}).
\end{aligned} \tag{5.21}
$$

The expansion of the second term of equation (5.17) of $\bar{R}^*$ is given by

$$
\begin{aligned}
&\tfrac{1}{6}(\hat{\psi} - \psi)^{-1} \hat{J}_{\psi\psi}^{-1/2} \sum_{r,s,t=1}^{k} \ell_{\theta_r\theta_s\theta_t}(\hat{\theta}) \hat{\delta}_r \hat{\delta}_s \hat{\delta}_t \\
&= -\tfrac{1}{6} V_{\psi\psi} Z_\psi^{-1} n^{1/2} (-V_{\psi\psi})^{-1/2} \\
&\quad \times n^{-1/2} (-1) V_{\psi\psi}^{-3} V_{\psi\psi\psi} Z_\psi^3 n^{-1/2} + O_p(1) \\
&= \tfrac{1}{6}(-V_{\psi\psi})^{-5/2} V_{\psi\psi\psi} Z_\psi^2 n^{-1/2} + O_p(n^{-1}).
\end{aligned} \tag{5.22}
$$

The expansion of the third term of equation (5.17) of $\bar{R}^*$ is given by

$$
\begin{aligned}
&\tfrac{1}{6}(\hat{\psi} - \psi)^{-3} \hat{J}_{\psi\psi}^{-3/2} \sum_{r,s,t=1}^{k} \ell_{\theta_r\theta_s\theta_t}(\hat{\theta}) \hat{\delta}_r \hat{\delta}_s \hat{\delta}_t \\
&= -\tfrac{1}{6}(-V_{\psi\psi})^{-3} Z_\psi^{-3} n^{3/2} (-V_{\psi\psi})^{-3/2} n^{-3/2} \\
&\quad \times (-V_{\psi\psi})^{-3} V_{\psi\psi\psi} Z_\psi^3 n^{-1/2} + O_p(n^{-1}) \\
&= -\tfrac{1}{6}(-V_{\psi\psi})^{-3/2} V_{\psi\psi\psi} n^{-1/2} + O_p(n^{-1}); \\
&= (-\ell_{\psi\psi}(\theta))^{-1} \ell_\psi(\theta) \ell_{\lambda\lambda\psi}(\theta) + O_p(1) \\
&= -V_{\psi\psi}^{-1} V_{\lambda\lambda\psi} Z_\psi \sqrt{n} + O_p(1).
\end{aligned} \tag{5.23}
$$

The expansion of the fourth term of equation (5.17) of $\bar{R}^*$ is given by

$$
\begin{aligned}
&-\tfrac{1}{2}(\hat{\psi} - \psi)^{-1} |\hat{J}_{\lambda\lambda}|^{-1} \hat{J}_{\psi\psi}^{-1/2} \operatorname{tr}\left( \operatorname{adj}(\hat{J}_{\lambda\lambda}) \sum_{r=1}^{k} \ell_{\lambda\lambda_r}(\theta) \hat{\delta}_r \right) \\
&= -\tfrac{1}{2}(-V_{\psi\psi}) Z_\psi^{-1} n^{1/2} (-V_{\psi\psi})^{-1/2} n^{-1/2} \\
&\quad \times \operatorname{tr}((-n^{-1}) V_{\lambda\lambda}^{-1} (-V_{\psi\psi}^{-1}) V_{\lambda\lambda\psi} Z_\psi \sqrt{n}) + O_p(n^{-1}) \\
&= \tfrac{1}{2}(-V_{\psi\psi})^{-1/2} \operatorname{tr}(V_{\lambda\lambda}^{-1} V_{\lambda\lambda\psi}) n^{-1/2} + O_p(n^{-1}).
\end{aligned} \tag{5.24}
$$

Finally, we have the following expansion of the fifth term of equation (5.17) of $\bar{R}^*$:

$$
\begin{aligned}
&\frac{\operatorname{tr}(\operatorname{adj}(A)X)}{R|A|}\\
&\quad= \frac{1}{R}\operatorname{tr}(A^{-1}X)\\
&\quad= -\frac{1}{2}(-\ell_{\psi\psi}(\theta))^{-3/2}\mu_{\psi,\psi\psi} + (-\ell_{\psi\psi}(\theta))^{-1/2}\operatorname{tr}(\ell_{\lambda\lambda}^{-1}(\theta)\mu_{\lambda,\lambda\psi})\\
&\qquad + O_p(n^{-1}) \qquad (5.25)\\
&\quad= -\frac{1}{2}n^{-3/2}(-V_{\psi\psi})^{-3/2}nv_{\psi,\psi\psi}\\
&\qquad + n^{-1/2}(-V_{\psi\psi})^{-1/2}\operatorname{tr}(V_{\lambda\lambda}^{-1}n^{-1}nv_{\lambda,\lambda\psi}) + O_p(n^{-1})\\
&\quad= \left[-\frac{1}{2}(-V_{\psi\psi})^{-3/2}v_{\psi,\psi\psi} + (-V_{\psi\psi})^{-1/2}\operatorname{tr}(V_{\lambda\lambda}^{-1}v_{\lambda,\lambda\psi})\right]n^{-1/2}\\
&\qquad + O_p(n^{-1}).
\end{aligned}
$$

Using the results in equations (5.21)–(5.25) in equation (5.17) yields

$$
\begin{aligned}
\bar{R}^* &= (-V_{\psi\psi})^{-1/2}Z_\psi\\
&\quad + [\tfrac{1}{2}(-V_{\psi\psi})^{-3/2}Z_{\psi\psi}Z_\psi\\
&\qquad - (-V_{\psi\psi})^{-1/2}Z_{\psi\lambda}^T V_{\lambda\lambda}^{-1}Z_\lambda\\
&\qquad - \tfrac{1}{2}(-v_{\psi\psi})^{-3/2}V_{\psi\psi\lambda}^T V_{\lambda\lambda}^{-1}Z_\lambda Z_\psi\\
&\qquad + \tfrac{1}{2}(-V_{\psi\psi})^{-1/2}Z_\lambda^T V_{\lambda\lambda}^{-1}V_{\psi\lambda\lambda}V_{\lambda\lambda}^{-1}Z_\lambda \qquad (5.26)\\
&\qquad + \tfrac{1}{6}(-V_{\psi\psi})^{-5/2}V_{\psi\psi\psi}Z_\psi^2 - \tfrac{1}{6}(-V_{\psi\psi})^{-3/2}V_{\psi\psi\psi}\\
&\qquad + \tfrac{1}{2}(-V_{\psi\psi})^{-1/2}\operatorname{tr}(v_{\lambda\lambda}^{-1}V_{\lambda\lambda\psi})\\
&\qquad - \tfrac{1}{2}(-v_{\psi\psi})^{-3/2}v_{\psi,\psi\psi}\\
&\qquad\qquad + (-V_{\psi\psi})^{-1/2}\operatorname{tr}(V_{\lambda\lambda}^{-1}v_{\lambda,\lambda\psi})]n^{-1/2} + O_p(n^{-1}),
\end{aligned}
$$

where

$$
v_{\psi,\psi\psi} = \frac{1}{n}E[\ell_\psi(\theta)\ell_{\psi\psi}(\theta);\theta], \qquad v_{\lambda,\lambda\psi} = \frac{1}{n}E[\ell_\lambda(\theta)\ell_{\lambda\psi}^T(\theta);\theta].
$$

This expansion has a remainder term of order $O_p(n^{-1})$. The same general approach can be used to obtain an expansion of $\bar{R}^*$ with a remainder of order $O_p(n^{-3/2})$. Such an expansion has the form

$$
\bar{R}^* = (-V_{\psi\psi})^{-1/2}Z_\psi
$$

$$
\begin{aligned}
&+ [\tfrac{1}{2}(-V_{\psi\psi})^{-3/2} Z_{\psi\psi} Z_{\psi} \\
&\quad - (-V_{\psi\psi})^{-1/2} Z_{\psi\lambda}^T V_{\lambda\lambda}^{-1} Z_{\lambda} \\
&\quad - \tfrac{1}{2}(-v_{\psi\psi})^{-3/2} V_{\psi\psi\lambda}^T V_{\lambda\lambda}^{-1} Z_{\lambda} Z_{\psi} \\
&\quad + \tfrac{1}{2}(-V_{\psi\psi})^{-1/2} Z_{\lambda}^T V_{\lambda\lambda}^{-1} V_{\psi\lambda\lambda} V_{\lambda\lambda}^{-1} Z_{\lambda} \\
&\quad + \tfrac{1}{6}(-V_{\psi\psi})^{-5/2} V_{\psi\psi\psi} Z_{\psi}^2 - \tfrac{1}{6}(-V_{\psi\psi})^{-3/2} V_{\psi\psi\psi} \\
&\quad + \tfrac{1}{2}(-V_{\psi\psi})^{-1/2} \operatorname{tr}(v_{\lambda\lambda}^{-1} V_{\lambda\lambda\psi}) \\
&\quad - \tfrac{1}{2}(-v_{\psi\psi})^{-3/2} v_{\psi,\psi\psi} \\
&\qquad\qquad + (-V_{\psi\psi})^{-1/2} \operatorname{tr}(V_{\lambda\lambda}^{-1} v_{\lambda,\lambda\psi})] n^{-1/2} \\
&+ R(n^{-1}) + O_p(n^{-3/2}),
\end{aligned} \tag{5.27}
$$

where $R(n^{-1})$ represents a finite number of terms of order $n^{-1}$. The specific form of $R(n^{-1})$ is omitted since it is extremely repetitive and tedious; it includes log-likelihood derivatives of the form $Z_{\theta_r\theta_s\theta_t}, r, s, t = 1, 2, \ldots, k$. Note that, since $R(n^{-1})$ represents a finite number of terms of order $n^{-1}$, it is a polynomial of $Z$ such that $E[R(n^{-1})] = O(n^{-1})$.

**Acknowledgments.** The authors would like to thank Editor Jianqing Fan, the Associate Editor and two anonymous referees for their very helpful and constructive suggestions and comments.

Department of Mathematics and Statistics
The University of Melbourne
Melbourne, Victoria 3010
Australia
E-mail: h.he@ms.unimelb.edu.au

Department of Statistics
Northwestern University
2006 Sheridan Road
Evanston, Illinois 60208
USA
E-mail: severini@northwestern.edu